\documentclass[11pt]{amsart}
\usepackage{amssymb}
\usepackage{amsthm}
\usepackage[T1]{fontenc}
\usepackage{enumerate}
\usepackage[colorlinks]{hyperref}
\usepackage{verbatim}
\usepackage{mathtools}
\usepackage{tikz-cd}
\usetikzlibrary{arrows}
\usepackage{algorithmic}
\usepackage[boxed]{algorithm}

\newcommand{\cN}{\mathcal{N}}

\newtheorem{theorem}{Theorem}[section]
\newtheorem{corollary}[theorem]{Corollary}

\newtheorem{lemma}[theorem]{Lemma}
\newtheorem*{TheoremNoN}{Theorem}

\theoremstyle{definition}
\newtheorem{definition}[theorem]{Definition}

\theoremstyle{remark}
\newtheorem{remark}[theorem]{Remark}
\newtheorem{example}[theorem]{Example}

\usepackage[latin1]{inputenc}
\usepackage[english]{babel}

\newcommand{\Pl}{Pl\"ucker}

\newcommand{\Proj}{\textnormal{Proj}\,}
\newcommand{\Supp}{\textnormal{Supp}\,}

\newcommand{\sat}{{\textnormal{sat}}}
\newcommand{\reg}{\textnormal{reg}}

\newcommand{\Ht}{\textnormal{Ht}}

\newcommand{\id}{\mathfrak}

\newcommand{\Mf}{\mathcal{M}\textnormal{f}}
\newcommand{\hilb}{{\mathcal{H}\textnormal{ilb}}}
\newcommand{\hilbp}{{\mathcal{H}\textnormal{ilb}_{p(t)}^n}}
 
\newcommand{\PS}{S} 
\newcommand{\PGL}{\textnormal{PGL}}
\newcommand{\GL}{\textnormal{GL}}
\newcommand{\PP}{\mathbb{P}} 
\newcommand{\Af}{\mathbb{A}}

\DeclareMathAlphabet{\mathpzc}{OT1}{pzc}{m}{it} 

\begin{document} 

\title{A Borel open cover of  the Hilbert scheme}

\author[C.~Bertone]{Cristina Bertone}
\address{Dipartimento di Matematica dell'Universit\`{a} di Torino\\ 
         Via Carlo Alberto 10, 
         10123 Torino, Italy}
\email{\href{mailto:cristina.bertone@unito.it}{cristina.bertone@unito.it}}

\author[P.~Lella]{Paolo Lella}
\email{\href{mailto:paolo.lella@unito.it}{paolo.lella@unito.it}}
\urladdr{\url{http://www.personalweb.unito.it/paolo.lella/}}

\author[M.~Roggero]{Margherita Roggero}
\email{\href{mailto:margherita.roggero@unito.it}{margherita.roggero@unito.it}}

\subjclass[2010]{14C05, 14Q20, 13P10} 
\keywords{Hilbert scheme, Borel-fixed ideal, marked scheme}

\begin{abstract} 
 Let $p(t)$ be an admissible Hilbert polynomial in  $\PP^n$ of degree $d$. The Hilbert scheme $\hilbp$ can be realized as a closed subscheme of a suitable  Grassmannian  $ \mathbb G$, hence  it could be globally defined by homogeneous equations in the \Pl\ coordinates of  $ \mathbb G$ and covered by open subsets given by the non-vanishing of a \Pl\  coordinate,  each embedded as a closed subscheme of the affine space $\Af^D$,   $D=\dim( \mathbb G)$.  However, the number $E$ of \Pl\ coordinates  is so large that  effective computations in this setting   are practically impossible.   In this paper, taking advantage of the symmetries of $\hilbp$, we exhibit a new open cover, consisting of marked schemes over  Borel-fixed ideals, whose number is significantly smaller than $E$. Exploiting the properties  of marked schemes, we  prove that these open subsets   are defined by equations of degree $\leq d+2$ in their natural embedding  in  $\Af^D$. Furthermore we find new embeddings in  affine spaces of  far lower dimension than $D$, and characterize  those that are still defined 
by equations of degree $\leq d+2$. The  proofs are constructive and use a polynomial reduction process,  similar to the one for Gr\"obner\  bases, but are term order free. In this new setting, we can achieve explicit computations in many non-trivial cases.
 \end{abstract} 

\maketitle

\section*{Introduction}

The Hilbert scheme, associated to the projective space
$\PP^n$ and to an admissible Hilbert polynomial $p(t)$ , was first introduced by  \cite{Gro} and parametrizes the set of
all subschemes $Z$ of $\PP^n$ with Hilbert polynomial $p(t)$. {We use $\hilbp$ for $\mathcal H\mathrm{ilb}_{p(t)}(\PP^n)$ throughout the paper.} 
The aim of this paper, and of some other related ones (see  \cite{CLMR,CR,LR}), is to find effective methods that allow  explicit computations on $\hilbp$. We choose the following definition of Hilbert scheme  since it is  the one that best suits  this purpose. 
  
Let $\PS = K[x_0,\dots, x_n]$ be a polynomial ring in $n + 1$ variables with coefficients in a field $K$ of characteristic $0$ and  let  $\PS_m$ be  the homogeneous component of $\PS$ of degree $m$. We denote by $r$ the Gotzmann number of $p(t)$, whose definition and properties were first established  by  \cite{gotz}, and by $q(t) $ the Gotzmann (or volume) polynomial  $N(t)- p(t)=\dim_K(\PS_t)- p(t)$. 
Every subscheme Z of $\PP^n$ with Hilbert polynomial $p(t)$ is  defined by a homogeneous saturated  ideal $\mathcal I_Z\subset \PS$  whose homogeneous component of degree $r$ is a $k$-vector space of dimension  $s:=q(r)$.  It can be proved that  the vector space $\mathcal I_Z\cap S_r$ generates $(\mathcal I_Z)_{\geq r}$ and then completely determines $Z$; by an abuse of notation we will write $I\in \hilbp$ meaning that $\Proj(S/I)\in \hilbp$. Further, using the correspondence $\mathcal I_Z \mapsto (\mathcal I_Z)_r$,   the Hilbert scheme can be realized as a closed subscheme of  the Grassmannian
 $ \mathbb G=\mathbb G(s,S_r)$ of linear spaces of dimension $s$ in $S_r$. A set of generators for the ideal defining $\hilbp$ in the ring of \Pl\ coordinates of $ \mathbb G$ can be obtained by imposing on an $s$-dimensional subspace  $V$ of $S_r$ the condition that the ideal $(V)\subset S$ has  the prescribed Hilbert polynomial $q(t)$. Explicitly,  by Gotzmann theorems of minimal growth \cite{gotz}, it is sufficient to impose that  $\dim_K S_1V$ is exactly $q(r+1)$. 
This setting dates back to \cite{Gro}, \cite{Mumford} (regularity) and \cite{gotz} (regularity and persistence), and it is fully described in \cite{IarroKanev} and also in \cite{CS} and  \cite{ACG}.

This way to  define the Hilbert scheme is constructive, because it relies on linear algebra. Nevertheless, explicit computations were not carried out until now, except in trivial cases, because the number $E$  of \Pl\ coordinates  is almost always huge. 
For instance, for  $n=3$ and $p(t)=4t$  the Gotzmann number is  $6$ and  $\hilb^3_{4t}$  becomes a subscheme of $\mathbb G(60,84)$; a rough  computation shows that   the number $E$ of \Pl\ coordinates is  $\sim 6\cdot 10^{20}$.

In  spite of these difficulties, many authors have dealt with the problem of determining a
set of explicit equations for $\hilbp$ as a subscheme of 
$ \mathbb G$.  For instance in \cite{IarroKanev},  \cite{B}, \cite{HaimSturm}, and more recently in \cite{ABM}  and \cite{BLMR} the authors determine upper bounds for the degree of  special sets of  equations in  the \Pl\ coordinates of $ \mathbb G$  defining scheme-theoretically   $\hilbp$.
In particular in \cite{BLMR} the authors obtain a set of equations for $\hilbp$  of  degree  $\leq \deg(p(t))+2$,   lower than the previous ones.  However, an explicit and computable presentation of $\hilbp$ is not yet easily achieved from this setting. \\
A simpler way to address this issue  is to consider an open cover of $\mathbb G$, consisting of  cells of maximal dimension $D=p(r)q(r)$,  and study the induced open cover of  $\hilbp$.
Each  \Pl\ coordinate can be naturally coupled to a monomial ideal $J\subseteq S$ generated by $s$ linearly independent monomials of degree $r$. If $\Delta_J$ denotes the  \Pl\ coordinate corresponding to $J$ and   $\mathcal U_J\simeq \Af^D$ is the open cell of  $\mathbb G$ given by $\Delta_J\neq 0$,  it will be sufficient to compute equations for $\mathcal H_J= \hilbp \cap \mathcal U_J$  as a closed subscheme of $\Af^D$. 

 However, in order to make this strategy truly effective, we must overcome two main difficulties:   the first one is finding an open cover of $\hilbp$ consisting of a  reasonably small number of open sets $\mathcal H_J$; the second one is finding an efficient way to compute equations for   $\mathcal H_J$. We are able to achieve both goals using   monomial ideals which have strong combinatorial properties and a long history in the theory of Hilbert schemes: the \emph{Borel-fixed ideals} (see Definition \ref{definBorel}). 

In Section \ref{sec:notation} we prove, for   monomial ideals $J$  of this type, the following
 result (Theorem \ref{isomHM}) that establishes a close relation  between  the open subset $\mathcal H_J$ of $\hilbp$ and the $J$-\emph{marked scheme}   $\Mf(J)$   introduced in \cite{CR}:

\begin{TheoremNoN} {Let   $J$ be  a Borel-fixed ideal, such that $S/J$ has Hilbert polynomial $p(t)$, let $r$ be the Gotzmann number of $p(t)$. Furthermore, assume that   $J={J^\sat}_{\geq r}$. Then   there is a commutative diagram of morphisms of  schemes:}
\begin{equation*}
\begin{tikzpicture}
\node (1) at (0,0) [] {$\mathcal H_{J}$};
\node (2) at (2,0) [] {$\Mf(J)$};
\node (3) at (0,-1.5) [] {$\mathcal U_J$};
\node (4) at (2,-1.5) [] {$\Af^D$};
\draw [->] (1) --node[above]{$\sim$} (2);
\draw [right hook->] (1) -- (3);
\draw [right hook->] (2) -- (4);
\draw [->] (3) --node[above]{$\sim$} (4);
\end{tikzpicture}
\end{equation*}
\end{TheoremNoN}
  The main results of the present paper  strongly rely 
 on the features and properties of marked schemes proved in \cite{CR} and \cite{BCLR}.

In Section \ref{sec:grass}, we show that the open subsets $\mathcal H_J $ corresponding to Borel-fixed ideals $J$ with Hilbert polynomial $q(t)$, are sufficient to define  an open cover of $\hilbp$, up to changes of coordinates in $\mathbb P^n$, while if $J$ is Borel but its Hilbert polynomial is not $q(t)$, then $\mathcal H_J=\emptyset$ (cf. Theorem \ref{fuori}). Furthermore we present an explicit algorithm that, given an ideal $I$ defining a point of $\hilbp$, returns a suitable change of coordinates $g$ of $\PP^n$ and a Borel-fixed ideal $J$ such that $I^g\in \mathcal H_J$. The algorithm exploits results of \cite{CLMR} and of \cite{L} to compute the list of Borel-fixed ideals belonging to $\hilbp$.

In Section \ref{quarta} we investigate the scheme structure of $\mathcal H_J$  based on the isomorphism with $\Mf(J)$ given in Theorem \ref{isomHM}.   The results about the theory of marked schemes obtained in  \cite {CR}  allow us to prove that   $\mathcal H_J$   can be  scheme theoretically  defined by equations of degree $\leq \deg(p(t))+2$ as a subscheme of $\Af^D$. The improvements introduced in \cite{BCLR} lead to the definition of   new embeddings of $\mathcal H_J$  in affine spaces of dimension  lower than  $D$ (cf. Theorem \ref{isomH}).   
Among them,  we determine those defined by equations of degree $\leq \deg(p(t))+2$ and compute upper bounds for the number of such defining equations (cf. Theorem \ref{gradod+2}); however the embedding corresponding to an affine space $\Af^{D'}$ of lowest dimension $D'$ may be given by equations of very high degree.

Summing up, we obtain an open cover of $\hilbp$ such that each open subset can be realized as a subscheme of a suitable affine space by a bounded number of  equations of degree $\leq \deg(p(t))+2$ (see Theorem \ref{gradod+2}) .

We emphasize that our proofs are constructive because they  are based on the features of some special  polynomial reduction processes introduced in \cite{CR} and in \cite{BCLR}. These are similar to ones for Gr\"obner\  bases, but are term order free: hence we can achieve explicit computations in many non-trivial cases. We apply our results to some examples in the final Section \ref{secexamples}.

\section{Hilbert scheme, Borel-fixed ideals and marked schemes} \label{sec:notation}

In this section we  first introduce the notations used throughout the paper and  present the three protagonists of the paper: first of all the Hilbert scheme parameterizing subschemes of $\mathbb P^n$ with a prescribed Hilbert polynomial, then the Borel-fixed ideals with their main features and finally the $J$-marked scheme, with $J$ a Borel ideal. At the end of the section (Theorem \ref{isomHM}) we will establish a first relation among them.

Throughout the paper $K$ is a field of characteristic 0, and 
 $\PS:=K[x_0,\dots,x_n]$ is the  polynomial ring over $K$  in the set of variables $x_0,\dots,x_n$. The elements and ideals that we consider  are always  homogeneous.

An ideal $J$ in $S$ is \emph{monomial} if it is generated by monomials. If $J\subset S$ is a monomial ideal, then it has a unique minimal set of monomial generators \cite[Lemma 1.2]{MS} that we call the \emph{monomial basis of $J$} and denote by $B_J$.

\subsection{\texorpdfstring{The Hilbert scheme $\hilbp$}{The Hilbert scheme}}\label{secHS}

In the following, we consider  a degree $d$ Hilbert polynomial $p(t)$  in the projective space  $\PP^n=\Proj(\PS)$. $\hilbp$ will denote the \emph{Hilbert scheme} parameterizing all subschemes $Z \subset \PP^n$ with Hilbert polynomial $p(t)$. We will denote by $r$ the Gotzmann number of $p(t)$, that is the highest Castelnuovo-Mumford regularity for the subschemes parameterized by $\hilbp$ (see \cite[Theorem 3.11]{Gr}), by $N(t)$ the binomial $\binom{n+t}{n}$ and  by $q(t)$ the volume polynomial (or Gotzmann polynomial) $N(t)-p(t)$. Moreover, we set 
   $s=q(r)$, $s'=q(r+1)$. 

\medskip

We  now recall one of the ways to embed the Hilbert scheme as a closed subscheme of a Grassmannian of suitable dimension. This setting dates back to \cite{Gro} and \cite{gotz} and it is fully described in  \cite{IarroKanev}; it can also be found in  \cite{CS} and  \cite[Chapter IX]{ACG}.   {
This is only one of the several ways to embed the Hilbert scheme $\hilbp$ in a Grassmannian $\mathbb G(p(r), S_r )$ \cite{Gro,IarroKanev,HaimSturm,ABM,BLMR} or product of Grassmannians $\mathbb G(p(r), S_r) \times \mathbb G(p(r + 1), S_{r+1})$ \cite{gotz, IarroKanev}. 
For a comparison between the different sets of equations defining such embeddings of $\hilbp$, see \cite{BLMR}.}

Every subscheme $Z\in \hilbp$ can be defined by many different ideals $I$ in $\PS$ such that $S/I$ has Hilbert polynomial $p(t)$; among them, there is the saturated ideal $\mathcal{I}_Z$, whose regularity is less than or equal $r$. We now consider $Z\in \hilbp$ as defined by    the truncated ideal $I:=(\mathcal{I}_Z)_{\geqslant r}$, which is generated by $s$  linearly independent forms of degree $r$, therefore it is uniquely determined by a linear subspace of dimension $s$ in the $K$-vector space $\PS_r$ of dimension    $N(r)$. Thus, $\hilbp$ can be set-theoretically embedded in the  Grassmannian $ \mathbb{G}:=\mathbb{G}(s,\PS_r)$ of all  $s$-dimensional subspaces of $\PS_r$.

By abuse of notation, we will write $I \in  \mathbb{G}$  if $I$ is the ideal generated by the vector space $I_r \in \mathbb G$ and we will write $I \in \hilbp$ if $I\in \mathbb{G}$ and its Hilbert polynomial is the volume polynomial $q(t)$. We will denote by $D$ the dimension of $\mathbb{G}$.

The scheme structure of $\hilbp$ is defined by imposing to $I\in \mathbb{G}$  that the dimension of the vector space $I_{r+1}$ is $s'$ (see \cite[Theorems C.17]{IarroKanev}). Moreover, by Macaulay's lower bound on growth,  the inequality  $\dim_K(I_{r+1})\geqslant s'$ is always true and therefore,   the condition $\dim_K(I_{r+1})= s'$   is in fact equivalent to  $\dim_K(I_{r+1})\leqslant s'$. 

{ If we fix  bases for the $K$-vector spaces $\PS_m$, $m\in \mathbb N$ (for instance the bases of monomials),  every vector space $V$ in $S_m$ can be represented by a (not unique) matrix $M(V,m)$ whose  rows contain the coefficients w.r.t. the fixed basis of $S_m$ of  a set of  polynomials that generate  $V$. In particular,   every ideal  $I\in \mathbb{G}$ can be represented by a  $s\times N(r)$ matrix $M(I_r,r)$, whose 
 minors of maximal order $s$  are the \Pl\ coordinates of $I$.  Moreover, by Macaulay's growth theorem, the rank of  the matrix $M(I_{r+1}=S_1 I_{r} ,r+1)$ in $S_{r+1}$    is  $\geq s'$. By the Gotzmann's Hilbert Scheme Theorem, equality holds if and only if $I\in \hilbp$. In this way  $\hilbp$ 
 can be defined by a homogeneous ideal in the ring of \Pl\ coordinates $K[\Delta]$ (see \cite{BLMR} for more details).
 
 In the following we fix the set of monomials as  basis of $\PS_m$ for every $m\in \mathbb N$.
 Every \Pl\ coordinate $\Delta$ of $I\in \mathbb{G}$ corresponds to a subset of $s$ elements of the fixed  basis of $S_r$}, so there is a one-to-one correspondence between the \Pl\ coordinates of the Grassmannian $\mathbb{G}$  and sets of $s$ monomials of degree $r$:  we will omit \lq\lq$r$\rq\rq\ and denote by $\Delta_J(I)$ (instead of $\Delta_{J_r}(I_r)$) the \Pl\ coordinate of $I\in \mathbb{G}$ corresponding to the monomial ideal $J\in \mathbb{G}$. In this paper \Pl\ coordinates corresponding to \emph{Borel-fixed monomial ideals} (introduced in the next section) will take center stage.

\subsection{Borel-fixed ideals}

We will use the compact notation $x$ for the set of variables $x_0,\dots, x_n$ and $x^\alpha$ for monomials in $\PS$, where $\alpha$ represents a multi-index $(\alpha_0,\dots,\alpha_n)$ of non-negative integers, that is $x^\alpha = x_0^{\alpha_0} \cdots x_n^{\alpha_n}$. 
The symbol $x^\alpha \mid x^\gamma$ means that $x^\alpha$ divides $x^\gamma$, that is there exists a monomial $x^\beta$ such that $x^\alpha \cdot x^\beta = x^\gamma$. If such monomial does not exist, we will write $x^\alpha \nmid x^\gamma$.

We consider the standard grading on the polynomial ring $\PS = \bigoplus_{m\in\mathbb{N}} \PS_m$, where $\PS_m$ denotes the homogeneous component of degree $m$; let $\PS_{\geqslant m} = \oplus_{m'\geqslant m}\left(\PS_{m'}\right)$ and in the same way, for every subset $\id{a}\subseteq \PS$, we will denote by $\id{a}_m$ and $\id{a}_{\geqslant m}$  the intersections of $\id{a}$ with $\PS_m$ and $\PS_{\geqslant m}$ respectively. 

We always order  the variables of $S$ in the following way: $x_0<x_1<\cdots<x_n$. 
\medskip

For every monomial $x^\alpha\not= 1$, we set $\min(x^\alpha)= \min\{x_i : x_i \mid x^\alpha\}$ and $\max(x^\alpha)= \max\{x_i : x_i \mid x^\alpha\}$.
We will say that a monomial $x^\beta$ can be obtained by a monomial $x^\alpha$ through an \emph{elementary  move} if   $x^\alpha x_j = x^\beta x_i$ for some variables $x_i\neq x_j$ or equivalently if there exists a monomial $x^\delta$ such that $x^\alpha=x^\delta x_i$ and $x^\beta=x^\delta x_j$. If $i < j$, we say that $x^\beta$ can be obtained by $x^\alpha$ through an \emph{increasing} elementary move and we write $x^\beta = \mathrm{e}_{i,j}^{+}(x^\alpha)$. The transitive closure of the  relation $x^\beta > x^\alpha$ if $x^\beta = \mathrm{e}_{i,j}^{+}(x^\alpha)$  gives a partial order on the set of monomials of any fixed degree (often called \emph{Borel partial order}), that we will denote by $>_{B}$:
\begin{equation*}
x^\beta  >_{B} x^\alpha\ \Longleftrightarrow\ \exists\  x^{\gamma_1}, \dots, x^{\gamma_t} \text{ such that } x^{\gamma_1} = \mathrm{e}^+_{i_0,j_0}(x^\alpha),\    \dots\ ,x^\beta = \mathrm{e}^{+}_{i_t,j_t}(x^{\gamma_t})
\end{equation*}
Note that every term order $\succ$ is a refinement of the Borel partial order $>_B$, that is $x^\beta >_B x^\alpha\ \Rightarrow\ x^\beta\succ x^\alpha$.

\begin{definition}\label{definBorel}
A monomial ideal $J \subset K[x_0, \dots, x_n]$ is said to be \emph{strongly stable} if every monomial $x^\alpha$ such that
$x^\alpha >_B x^\beta$, with $x^\beta \in J$, belongs to $J$.
\end{definition}

Under the hypothesis that $\textnormal{char}(K)=0$, a monomial ideal $J$ is strongly stable if, and only if,  $J$ is fixed by the action of the Borel subgroup of $\mathrm{GL}(n+1)$  of lower triangular matrices (see for instance \cite[Proposition 1.25]{Gr} for a proof). Using this equivalence, we apply to the \emph{Borel-fixed ideals} or, for short,  \emph{Borel ideals},  results proved  in \cite{CR} and \cite{BCLR} for strongly stable ideals.
{ We consider the usual action of $\mathrm{GL}(n+1)$ on the ring $S$ given by the change of coordinates $x_i \mapsto \sum_j g_{ij} x_j$ for any invertible matrix $g = (g_{ij})$. For any polynomial $f(x_0,\ldots,x_n)$ and any change of coordinates $g$, we denote by $f^{g}$ the result of the action of $g$ over $f$, i.e. $f(\ldots,\sum_j g_{ij} x_j,\ldots)$. 
In the same way, given an ideal $I \subset S$, we denote by $I^g$ the ideal $(f^g\ \vert\ f \in I)$.  }

The main reason why we assume that $K$ has characteristic 0, is Galligo's Theorem, which is a key point for our investigations (mainly in the proof of Lemma \ref{lem:ricgrass}).
\begin{theorem}\label{Galligo}{(\cite{Gall})}
Assume that $\mathrm{char}\,(K)=0$ and fix a term order $\prec$ on the monomials of $\PS$. For every ideal $I\subset \PS$ there exists a Zariski open subset $U\subset \mathrm{GL}(n+1)$ such that for every $g\in U$, the initial ideal $\mathrm{in}_\prec(I^g)$ is constant and strongly stable.
\end{theorem}

\begin{remark}\label{rational}
Since $K$ contains the rationals $\mathbb Q$, and $\GL(n+1,\mathbb Q)$ is Zariski-dense in  $\GL(n+1,K)$, we may restrict $g$ to a general change of coordinates with coefficients in $\mathbb Q$. This can be useful in practical applications, as in direct computation the use of rational numbers is more efficient.
\end{remark}

Thanks to  Galligo's Theorem,   any component, and any intersection of components on the Hilbert scheme, contains a point corresponding to a Borel ideal. Therefore, it   is   natural  to investigate the structure of $\hilbp$  using the Borel ideals and the action of the linear group. This will be our approach and  was  recently considered for instance by \cite{Sh}  to prove a first order infinitesimal version of Galligo's Theorem.

A homogeneous ideal $I$ is $m$-{\it regular} if the $i$-th syzygy module of $I$ is generated in degree $\leq m+i$, for all $i\geq 0$. The {\it regularity} $\reg(I)$ of $I$
is the smallest integer $m$ for which $I$ is $m$-regular. The {\it saturation} of a homogeneous ideal $I$ is $I^{\sat}=\{f\in S\ \vert \ \ \forall \ j=0,\ldots,n,\exists \ r \in {\mathbb N} : x_j^r f \in I\}$. The ideal $I$ is {\it saturated} if $I^{\sat}=I$. 

We recall that if $J$ is Borel then $\reg(J) = \max \{\deg x^\alpha \ \vert\ x^\alpha \in B_J\}$ \cite[Proposition 2.9]{BS} and its saturation $J^\sat$ is $(J:x_0^\infty)$ (for example, see \cite[Corollary 2.10]{Gr}). Hence $\reg(J^\sat)\leq \reg(J)$. Furthermore if $B_J$ is the monomial basis of a Borel ideal $J$, then the saturation of $J$ is generated by $B_J\vert_{x_0=1}$, that is one deletes the variable $x_0$ in each monomial of $B_J$ \cite[Theorem 3.2]{Sh}.

For a monomial ideal $J\subset \PS$, we denote by $\cN(J)$ the \emph{sous-escalier} of the monomial $J$, that is the set of monomials not belonging to $J$.

\begin{lemma}\label{potenze}  Let $J$ be a Borel ideal in $S$. Then:
\begin{enumerate}[(i)]
\item $x^\alpha \in J\setminus B_J\ \Rightarrow\ \dfrac{x^\alpha}{\min(x^\alpha)} \in J$;
\item \label{potenze_ii}$x^\beta \in \cN(J)$ and  $x_ix^\beta \in J \ \Rightarrow$ either $x_ix^\beta \in B_J$ or $x_i > \min(x^\beta)$.
\item \label{potenze_iii} if  $J\in \hilbp$, then
$K[x_{0}, \dots, x_{d}]\subset \cN(J^\sat)$ and  $K[x_{d+1}, \dots, x_{n}]_{\geq m}\subset J^\sat$   for some $m\leq \reg(J^\sat)$.
\end{enumerate}
\end{lemma}

\begin{proof}\ 
\begin{enumerate}[(i)]
\item We consider $x^\alpha\in J\setminus B_J$. Then $x^\alpha=x^\gamma\cdot x_i=x^\delta\cdot x_j$, with $x^\gamma\in J$, $x_j=\min(x^\alpha)$. Then $x^\delta=e^+_{j,i}(x^\gamma)$,  and then $x^\delta\in J$, since $J$ is Borel. 
\item If $x^\beta$ belongs to $\cN(J)$,   and $x_ix^\beta\in J\setminus B_J$, then $x_i>\min(x^\beta)$ by the previous item.
\item 
If $x_i$ is the maximal variable that is not nilpotent in $\PS/J^\sat$, then  $K[x_{0}, \dots, x_i]\subset \cN(J^\sat)$ by the strongly stable property of  $J^\sat$. Hence $K[x_{0}, \dots, x_i] \hookrightarrow \PS/J^\sat$ so that $i \leqslant \dim{(\PS/J)} = d$. \\
On the other hand,  some power $x_{i+1}^m$ belongs to $B_J$; hence, again by the strongly stable property,  $K[x_{i+1},\dots,x_n]_{\geq m}\subset J$, so that   $i+1 > \dim (\PS/J^\sat) = d$.
\end{enumerate}
\end{proof}

\subsection{Marked schemes over Borel ideals}

In this subsection, we introduce the main definitions and properties concerning marked sets of polynomials and marked bases. These special sets were investigated in \cite{CR} and \cite{BCLR}, highlighting their interesting features under the hypothesis they are marked over a Borel ideal.

\begin{definition}
For any non-zero  polynomial $f \in \PS$, the \textit{support} of $f$ is the set of monomials $\Supp(f)$  that appear in $f$ with a non zero coefficient. By definition $\Supp(0)=\emptyset$.\\
A \emph{marked polynomial} is a polynomial $f\in \PS$ together with a specified monomial of $\Supp(f)$ that will be called \emph{head term} of $f$ and denoted by $\Ht(f)$. \\
Given a monomial ideal $J$ and an ideal $I$ such that $\cN(J)$ generates $S/I$ as a $K$-vector space, a {\em $J$-reduced form modulo $I$} of a  polynomial $h$ is a polynomial $h_0$ such that $h-h_0\in I$ and $\Supp(h_0) \subseteq \cN(J)$. If the $J$-reduced  form modulo $I$ is unique, we call it \emph{$J$-normal form modulo $I$}.
\end{definition}

Note that every polynomial $h$ has a unique $J$-reduced form modulo  $I$ if, and only if, $\cN(J)$ is a $K$-basis for the quotient $S/I$. In this case,  the $J$-reduced form modulo the homogeneous ideal $I$ of a homogeneous polynomial $h$ turns out to be  homogeneous too of the same degree. If we only assume that $\cN(J)$ generates $\PS/I$,  there could be  several  $J$-reduced forms modulo $I$ of  $h$, but among them  we can always find at least one  that is homogeneous.  In fact, if $h \in S_m$ and $h'$ is a $J$-reduced form modulo $I$ of $h$, then $(h-h')_m=h-h'_m$ belongs to $I$ too and  $\Supp(h'_m)\subseteq \Supp(h')\subseteq \cN(J)$. 
 Then $h'_m$ is a homogeneous $J$-reduced form modulo $I$ of $h$. Thus, we can always consider homogeneous  $J$-reduced forms of  homogeneous polynomials.

\begin{definition}\label{schemestruct}
A finite set $G$ of homogeneous marked polynomials $f_\alpha=x^\alpha-\sum c_{\alpha\gamma} x^\gamma$, with  $\Ht(f_\alpha)=x^\alpha$, is called $J$-\emph{marked set} if the head terms $\Ht(f_\alpha)$ are pairwise different, they form the monomial basis $B_J$ of the monomial ideal $J$ and every $x^\gamma$ belongs to $\cN(J)$, i.e.  $\vert\Supp(f)\cap J \vert =1$. We call \emph{tail} of $f_\alpha$ the polynomial $T(f_\alpha):=\Ht(f_\alpha)-f_\alpha$, so that $\Supp(T(f_\alpha))\subseteq \mathcal N(J)$.  A $J$-marked set $G$ is a $J$-\emph{marked basis} if $\cN(J)$ is a basis of $S/(G)$ as a $K$-vector space. 

The family of all homogeneous ideals $I$ such that $\cN(J)$ is a basis of the quotient $S/I$ as a $K$-vector space will be denoted by $\Mf(J)$. 
\end{definition}

If $J$ is a Borel ideal, then $\Mf(J)$ can be endowed with a natural structure of affine scheme (shown in \cite[Section 4]{CR} and recalled in Remark \ref{struttMf} below) called $J$-\emph{marked scheme}, that can be explicitly computed, using   a polynomial reduction process,  similar to the one for Gr\"obner\  bases, but term order free. Observe that  $\Mf(J)$   contains every homogeneous ideal having $J$ as initial ideal with respect to some term order, but in general it can also contain other ideals \cite[see][Example 3.18]{CR}.

As a  straightforward consequence of the  definition, the ideals of $\Mf(J)$ define points on the Hilbert scheme $\hilbp$ containing $\PS/J$. However, two different ideals  of 
 $\Mf(J)$ may correspond to the same  point \cite[Example 3.4]{BCLR}, while we get a one-to-one correspondence between $\Mf(J)$ and a subset of $\hilbp$ assuming that   $J$ is  an  \emph{$m$-truncation ideal} of its saturation. \cite[Theorem 3.3]{BCLR}.

\begin{definition} Let $I\subseteq S$ be a homogeneous  ideal. We will say that $I$ is an \emph{$m$-truncation} if $I$ is the truncation of $ I^{\sat}$ in degree $m$, that is  $I={I^\sat}_{\geq m}$.
\end{definition} 

\begin{remark}\label{ossTronc}
Observe that not every homogeneous  ideal $I$ generated by its degree $m$ component is an $m$-truncation. For instance the ideal $I=(x_0^2, \dots, x_i^2, \dots, x_n^2)$ is not a $2$-truncation, since $I^{\sat}=(1)$ and $I\neq \PS_2$. Nevertheless, in our setting of Section  \ref{secHS} every  ideal  $I\in \hilbp$ is an $r$-truncation. In fact, of course $I\subseteq (I^{\sat})_{\geq r}$. Moreover   the two ideals $I$ and $I^{\sat}$ have the same Hilbert polynomial $p(t)$ and their Hilbert functions coincide with $p(t)$ for every $t\geq r\geq \reg(I^{\sat})$.

Moreover, every   Borel  ideal $J\in \mathbb G$, as in Section \ref{secHS},  is an $r$-truncation. Indeed, $J^{\sat}=(B_J\vert_{x_0=1})$ and then  for every  $x^\gamma \in B_{J^\sat}$, the monomial  $x^\gamma x_0^{r-\vert \gamma\vert }\in B_J$. If $x^\beta$ belongs to ${J^\sat}_{\geq r}$, then $x^\beta=x^\gamma x^\delta$ for some $x^\gamma\in B_{J^\sat}$. Hence, $x^\beta>_B x^\gamma x_0^m$, $m\geq r-\vert \gamma\vert$,  and $x^\beta$ belongs to $J$ by the Borel property.

\end{remark}

\subsection{\texorpdfstring{Marked schemes as open subsets of $\hilbp$}{Marked schemes as open subsets of the Hilbert scheme}}

In the present subsection, we will consider a Borel ideal $J\in  \hilbp$  and show that the corresponding $J$-marked scheme is scheme theoretically isomorphic to an open subset of $\hilbp$.

We denote by $\mathcal{Q}$ the set of monomial ideals in $\mathbb{G}=\mathbb{G}(s,\PS_r)$ and  by $\mathcal{B}$ the set of  those that are Borel.  Moreover  for every $J\in \mathcal Q$,  $\mathcal{U}_{J}$ will be the open subset of $\mathbb{G}$ containing the ideals $I$ such that their \Pl\ coordinate $\Delta_J(I)$ is not $0$, and $\mathcal{H}_{J}:=\mathcal{U}_{J}\cap \hilbp$ will be the corresponding open subset of $\hilbp$.
 As is well known, $\mathcal{U}_J$ is a  Schubert cell of maximal dimension $D:={p(r)q(r)}=\dim (\mathbb{G})$ and is  isomorphic to the affine space $\mathbb{A}^D$. 

\begin{lemma}\label{lem:pluckerlocali} Let $I$ be an ideal in $\mathbb{G}$,  and  let $B_J$ be the  monomial basis of an ideal   $J\in \mathcal{Q}$. Then the following are equivalent:
\begin{enumerate}[(i)]
\item\label{it:pluckerlocali_i} $\Delta_J(I) \neq 0$:
\item\label{it:pluckerlocali_ii} $I_r$ can be represented by a matrix of the form $\left(\begin{array}{l|r} \mathrm{Id}  & \mathcal {C}\end{array}\right)$, where  the left block is the $s\times s$ identity matrix and corresponds to the monomials in $B_J$ and the entries of the right block $\mathcal{C}$ are constants $-c_{\alpha \gamma}\in K$, where $x^\alpha \in B_J$ and $x^\gamma \in \mathcal{N}(J)_r$;
\item\label{it:pluckerlocali_iii} $I$ is generated by a $J$-marked set:
\begin{equation*}
G=\left\{f_\alpha=x^\alpha-\sum c_{\alpha\gamma} x^\gamma : Ht(f_\alpha)=x^\alpha\in B_J\right\}
\end{equation*}
\end{enumerate}
\end{lemma}
\begin{proof}\ 

 {(\ref{it:pluckerlocali_i})} $\Rightarrow$ {(\ref{it:pluckerlocali_ii})}. Up to rearranging the columns, it is sufficient to multiply any matrix $M(I,r)$ (as in Section \ref{secHS}) by the inverse of its submatrix made up of the columns  corresponding to $B_J$, since its determinant is $\Delta_J(I)\neq 0$.

{(\ref{it:pluckerlocali_ii})} $\Rightarrow$ {(\ref{it:pluckerlocali_i})} is obvious. 

{(\ref{it:pluckerlocali_ii})} $\Leftrightarrow$ {(\ref{it:pluckerlocali_iii})}. The generators of $I$ given by the rows of $M(I,r)$ as in (\ref{it:pluckerlocali_ii}) are indeed a $J$-marked set and, conversely, the matrix containing the coefficients of the polynomials $f_{\alpha}$ has precisely the shape required in {(\ref{it:pluckerlocali_ii})}.
\end{proof}

\begin{definition}\label{defin_A}
Let $J$ be a Borel ideal, let $C = \{C_{\alpha\gamma}\ \vert\ x^\alpha \in B_J,\ x^\gamma \in \mathcal{N}(J)_{\vert\alpha\vert}\}$ denote a set of new variables and consider the ring $K[C,x]$. We denote by $\mathcal G \subset K[C,x]$ the $J$-marked set:
\begin{equation}\label{JbaseC} \mathcal{G}= \left\{F_\alpha=x^\alpha-\sum C_{\alpha\gamma} x^\gamma : Ht(F_\alpha)=x^\alpha\in B_J, x^\gamma \in \cN(J)\right\}\end{equation}
and by $\id{I}_J$ the ideal  generated by $\mathcal G$ in the ring  $K[C,x]$.
\end{definition}

\begin{corollary}\label{cor:newVariables} 
 Under the hypotheses of Lemma \ref{lem:pluckerlocali}, there is an isomorphism  between $\mathcal{U}_J$ and $\mathbb{A}^D$ $=\mathrm{Spec}\, (K[C])$ such that
the (closed) points in    $ \mathcal{U}_J$ correspond to all ideals that we obtain from $\id{I}_J$ by specializing the variables $C_{\alpha \gamma}$ to $c_{\alpha \gamma}\in K$.
\end{corollary}
\begin{proof}
Under the hypotheses of Lemma \ref{lem:pluckerlocali}, it is sufficient to fix an isomorphism $\mathcal{U}_J\cong \mathbb{A}^{D}$ such that the constants $c_{\alpha \gamma}$ are the coordinates of $I$ in $\mathbb{A}^{D}$. 
\end{proof}

\begin{remark}\label{struttMf}
$\Mf(J)$ is naturally endowed of the structure of an affine subscheme of $\mathbb A^{D}$. More specifically,  as shown by \cite[Lemma 4.2]{CR}, we can obtain a set of generators for the ideal $\mathfrak A_J$  in $K[C]$  defining $\Mf(J)$ by imposing  conditions on the rank of some matrices. We obtain the same structure by a polynomial reduction process similar to the Gr\"obner one, but term order free \cite[Theorem 3.12, Appendix]{CR}.
  \end{remark}

We now fix  $J\in \mathcal B\cap \hilbp$ and    show that $\mathcal H_J$ is nothing but  the $J$-marked scheme $\Mf(J)$.

\begin{theorem}\label{isomHM}
For $J\in \mathcal B\cap \hilbp$, there are scheme theoretic isomorphisms  $\phi_1$ and $\phi_2$  and embeddings such that the following diagram commutes:
\begin{equation}\label{diagramma}
\begin{tikzpicture}
\node (1) at (0,0) [] {$\mathcal H_{J}$};
\node (2) at (2,0) [] {$\Mf(J)$};
\node (3) at (0,-1.5) [] {$\mathcal U_J$};
\node (4) at (2,-1.5) [] {$\Af^{D}$};
\draw [->] (1) --node[above]{\scriptsize $\phi_1$} (2);
\draw [right hook->] (1) -- (3);
\draw [right hook->] (2) -- (4);
\draw [->] (3) --node[above]{\scriptsize $\phi_2$} (4);
\end{tikzpicture}
\end{equation}
 \end{theorem}  
 \begin{proof}
By Gotzmann's Hilbert Scheme Theorem, we can give to the open subset $\mathcal H_J$ the structure of closed affine subscheme of $\mathcal U_J\simeq \mathbb A^{D}$. Indeed, $I\in \mathcal U_J$ belongs to $\mathcal H_J$ if $\dim I_{r+1}=q(r+1)$. This can be obtained by considering the matrix of the coefficients of the generators of $I_{r+1}$ and imposing that its rank is $q(r+1)=\dim J_{r+1}$. This gives the scheme-theoretical embedding of $\mathcal H_{J}$ in $\mathcal U_J$ of the diagram (see also \cite[Proposition C.30]{IarroKanev}).

As recalled by Remark \ref{struttMf},  $\Mf(J)$ can be considered as an affine subscheme of $\mathbb A^{D}$: this gives the  second embedding of the diagram.   In the present  hypothesis,  we obtain a set of generators for the ideal $\mathfrak A_J$  in $K[C]$  defining $\Mf(J)$ by imposing  that the rank of the matrix corresponding to the degree $r+1$ is  $\leq \dim J_{r+1}$. This matrix turns out to be indeed  $M(\id I_J,{r+1})$ and $\dim J_{r+1}=q(r+1)$. This gives the isomorphism $\phi_1$ between $\mathcal H_J$ and $\Mf(J)$ of the diagram. 

Finally the isomorphism $\mathcal U_J\simeq \mathbb A^{D}$ is the one of Corollary \ref{cor:newVariables}.
 \end{proof}  
 
In the next section, we will further investigate  the open subsets $\mathcal H_J$ of $\hilbp$  for a Borel ideal $J\in \hilbp$, showing that they are sufficient to give a cover  of $\hilbp$ (up to the action of the group $\mathrm{PGL}(n+1)$). In Section \ref{quarta} we will use the tools and results of \cite{BCLR} to obtain embeddings of  $\mathcal H_J$ in affine spaces of dimension smaller than $D=p(r)q(r)$ and  find  upper bounds for the number and  degrees of the corresponding   equations.

\section{The Borel cover}\label{sec:grass}

In this section we will   investigate the families of open subsets $\mathcal U_J\subseteq \mathbb G$  and $\mathcal H_J\subseteq \hilbp$ and deduce from them a open cover of $\hilbp$. For a different approach giving a cover by locally closed subschemes of the Hilbert scheme of points see \cite{Lederer}.

It is quite obvious that they cover respectively $\mathbb G$ and $\hilbp$ as $J$ varies in $\mathcal Q$. However $\mathcal Q$  in general has a large number of elements.  In the study of $\hilbp$, it is quite natural to consider Borel ideals, since every component and any intersection of components of $\hilbp$ contains a point corresponding to a Borel ideal. Further
 $\mathcal B$ is a  comparatively small subset of $\mathcal Q$; then it would be convenient to find a way to cover $\mathbb G$ and $\hilbp$  by  subfamilies  indexed by $\mathcal B$.

The open subsets $\mathcal U_J$, as $J$ varies in $\mathcal B$,  do not cover the whole Grassmannian: for instance, any monomial ideal $J'\in \mathcal Q\setminus \mathcal B$ does not belong to $\mathcal U_J$, for every $J\in \mathcal B$. However,    thanks to Galligo's Theorem (Theorem \ref{Galligo}), we get a complete cover of $\mathbb{G}$ up to the action on  $\mathbb G$ induced by  the projective linear group $\PGL(n+1)$ ($\PGL$ for short) of changes of coordinates in $\PP^n$. As usual we denote by a superscript $g$  the action of an element $g$ of the group.

\begin{lemma}\label{lem:ricgrass} The open subsets  $\mathcal{U}_{J}$, $J\in \mathcal B$, up to the action of $\PGL$,  cover  $\mathbb{G}$, that is:
\begin{equation*}
\bigcup_{\begin{subarray}{c}J\in \mathcal{B}\\g\in \PGL\end{subarray}} \mathcal{U}_{J}^g=\mathbb{G}.
\end{equation*}
\end{lemma}
\begin{proof} 
Let $I \in \mathbb{G}$ be any ideal and let $\preceq$ be any term order on the monomials of $\PS$. Due to Galligo's Theorem \ref{Galligo}, in generic coordinates the  initial ideal $J'$  of $I$  w.r.t $\prec$ is Borel, and then $J:=(J'_r)$ is Borel too. Moreover, by construction  $J$ is generated by $s$ monomials of degree $r$ and then $J\in \mathcal{B}$. Hence  for a general $g\in \PGL$ we have $\Delta_J(I^g)\neq 0$ that is   $I^g\in \mathcal{U}_{J}$, so that $I\in \mathcal{U}^{g^{-1}}_{ J}$.   
\end{proof}

\begin{remark}
In the proof of Lemma \ref{lem:ricgrass} we deal with the generic initial ideal $J'$,  that may also have minimal generators of degree $>r$, so that sometimes  $J'\notin \mathbb{G}$. To avoid this problem, we consider instead the ideal $J$ generated by $J'_r$, which is Borel  and certainly belongs  to $\mathbb{G}$ because a basis for  $J'_r$ is given by   $s$ monomials in degree $r$. The same expedient is used in \cite{BLMR}.

We can rephrase Lemma \ref{lem:ricgrass} saying that for every $I\in \mathbb G$ there are: a Borel ideal  $J\in \mathcal{B}$ and  a generic linear change of coordinates $g\in \PGL$, such that  $I^g$ is generated by a $J$-marked set $G$: for instance we can choose $J=(\mathrm{in}_\prec(I^g)_r)$, for any fixed term order $\prec$ and a general $g$. 
Nevertheless,   $G$ does not need to be a  $J$-marked basis for $I^g$, namely  $I$ and $J$ may have different  Hilbert polynomials.
\end{remark}

\begin{example}\label{diffUH}
Let us consider the ideal $I=(x_2^2,x_1^2)$ in $K[x_0,x_1,x_2]$ belonging to $\mathbb G(2,6)$. After the change of coordinates  $g$: $x_2\rightarrow x_2$, $x_1\rightarrow x_2+x_1$, $x_0\rightarrow x_0$, we get $I^g=(x_2^2,x_2^2+2x_2x_1+x_1^2)$ whose initial ideal with respect to any term order such that $x_2>x_1>x_0$ is $J'=(x_2^2, x_2x_1, x_1^3)$. In fact,  $K$  being  a field of characteristic 0, the reduced Grobner basis  of $I^g$ is the set of three polynomials $\{f_1,f_2,f_3\}$, where  $f_1=x_2^2$, $f_2=x_2x_1+x_1^2/2$ and $f_3=x_1^3=4x_1 f_1-(4x_2-2x_1)f_2$.  Obviously $J'$ does not belong to $\mathbb G(2,6)$.
 Following the line of the proof of  Lemma \ref{lem:ricgrass} we then consider $J=(J'_2)=(x_2^2,x_2x_1)$: now $J\in \mathbb G(2,6)$ and moreover $I^g\in \mathcal U_J$, as one can see by the coefficient matrix of $\{f_1,f_2\}$, \eqref{matrice},  even though the Hilbert polynomial of $S/I$  is $4$ and that of $S/J$ is $t+2$.
\begin{equation}\label{matrice}
\left(\begin{array}{cc|cccc}
1& 0&0&0&0&0\\
0&1&1/2&0&0&0
\end{array}\right)
\end{equation}
\end{example}

We now restrict the above open cover to  the the Hilbert scheme $\hilbp$ and  denote by ${\mathcal B}_{p(t)}$ the Borel ideals $J$ such that $S/J$ has Hilbert polynomial $p(t)$. Since  $\hilbp \hookrightarrow \mathbb{G}$ is  a closed embedding,  the intersection $\mathcal{H}_J=\mathcal{U}_J\cap \hilbp$ is an open subset of $\hilbp$, which may possibly be empty. If $J\in \mathcal B_{p(t)}$, then of course $\mathcal{H}_J$ is non-empty because it contains $J$ itself.
We show now that the converse is also true and this justifies the following:

\begin{definition} The \emph{Borel cover} of $\hilbp$ is the family of all the open subsets  $\mathcal{H}^g_J$ with  $J\in \mathcal{B}_{p(t)}$ and $g\in  \PGL$. 
\end{definition}

\begin{theorem}\label{fuori} 
 If $J\in \mathcal{B}$, then:
\begin{equation*}
 \mathcal{H}_J\neq  \emptyset\quad\Longleftrightarrow\quad {J\in \mathcal B_{p(t)}}.
\end{equation*}
As a consequence, we get:
\begin{equation*}
\hilbp=\bigcup_{\begin{subarray}{c}g\in \PGL \\ J\in \mathcal{B}_{p(t)}\end{subarray}} \mathcal{H}^g_J.
\end{equation*} 
\end{theorem} 
\begin{proof} 
We prove only the non-trivial part $(\Rightarrow)$ of the first statement. Assume that $J\notin \hilbp$. By Gotzmann's Hilbert Scheme Theorem, this is equivalent to $\dim_K(J_{r+1})>q(r+1)$. If $I$ is any ideal in $\mathcal{U}_J$, then it has a set of generators as those given in Lemma \ref{lem:pluckerlocali} (\ref{it:pluckerlocali_iii}), so that  $\dim_K(I_{r+1})\geqslant \dim_K(J_{r+1})>q(r+1)$ (see \cite[Corollary 2.3]{CR}). Hence $I \notin \hilbp$.

The other  statement is a direct consequence of the first one and of Lemma \ref{lem:ricgrass}.
\end{proof}

\begin{corollary}
If $J$ belongs to $\mathcal B$, then: 
\begin{itemize}
\item the ideal $I$ belongs to $\mathcal U_J$ if, and only if, $I$ is generated by a  $J$-marked set;
\item the ideal $I$ belongs to $\mathcal H_J$ if, and only if,  $J\in \mathcal B_{p(t)}$ and $I$ is generated by a  $J$-marked basis.
\end{itemize}
\end{corollary}

 There are several strategies that we can follow in order to compute for a given ideal  $I\in \hilbp$ a change of coordinates  $g \in \PGL$ and a Borel ideal  $J\in \mathcal B_{p(t)}$ such that $I_r^g\in \mathcal H_{J}$. 

First of all, as observed in Remark \ref{rational}, we can consider a generic change of coordinates with rational coefficients.

The method used in the proof of Lemma \ref{lem:ricgrass} may be used to determine $J\in \mathcal B_{p(t)}$ . Note that we do not need to perform the Buchberger's Algorithm on the generators of $I^g$ in order to derive a complete Gr\"obner basis of $I^g$ , but we can simply compute, by means of a Gaussian reduction, $s=q(r)$ distinct initial terms with respect to any term order $\preceq$
of degree $r$ forms in $I^g$.  The ideal $J$ generated by  these monomials  is indeed the initial ideal $J'$ of $I^g$ with respect to $\preceq$.
 In fact, using Lemma  \ref{lem:pluckerlocali}, $I^g$ is generated by a $J$-marked set, hence $\dim J_{r+1} \leq \dim {I^g}_{r+1} = q(r + 1)$, by \cite[Corollary 2.3]{CR}. On the other hand, by Macaulay's growth theorem, $\dim J_{r+1}\geq q(r+1)$. {Then equality holds and by Gotzmann's Hilbert Scheme Theorem $J$ has the same Hilbert polynomial $q(t)$ as $I^g$ and ${J'}$. 
By construction  $J=(J_r)\subseteq  J'$; furthermore, for all $t< r$, we have that $\dim J_t=\dim {J'}_t=0$, while for all $t\geq r$, $\dim J_t=\dim {J'}_t=q(t)$; thus $J = {J'}$.}

Another strategy to determine a Borel ideal $J$, such that $\mathcal{H}_J$ contains the point $I^g$ (avoiding to compute the generic initial ideal), is that of computing the \Pl\ coordinate $\Delta_J(I)$ for every $J\in\mathcal B_{p(t)}$: given $n$ and $p(t)$, the complete list $\mathcal B_{p(t)}$ can be explicitly computed by the algorithms presented by \cite{CLMR} or the ones presented by \cite{L}. Starting from   this  list,  the following algorithm returns a Borel ideal $J$ in $\mathcal B_{p(t)}$ with minimum regularity such that $I^g\in \mathcal H_J$. The computational advantage of the condition about  the regularity will be discussed in Section \ref{quarta}. Further, an easy variant of the algorithm  finds $g\in \PGL$ and the maximal subset  $\mathcal B(I^g)\subseteq \mathcal B_{p(t)}$ such that $I^g\in\mathcal H_{J}$ for every $J\in \mathcal B(I^g)$.

Let us suppose that the following functions are made available.
\begin{itemize}
\item  $\textsc{SuitableGrasmannian}(I)$. It computes the suitable Grassmannian (in the sense of Section \ref{secHS}) in which the Hilbert scheme $\hilbp$ is embedded, where $p(t)$ is the Hilbert polynomial of $S/I$. It returns the pair $(\tilde I, \mathbb G)$ such that $\tilde I$ is the point of $\mathbb G$ corresponding to $I$.
\item $\textsc{QGenericChangeOfCoordinates}(I)$. It performs a random linear change of coordinates $g$ with rational coefficients on the generators of the ideal $I$, returning the pair $(g,I^g)$.
\item $\textsc{BorelIdeals}(n,p(t))$. It returns the full list of Borel ideals $J$ in $\mathcal B_{p(t)}$; this list is sorted according to encreasing regularity of $J^\sat$.
\item $\textsc{Pl\"uckerCoordinate}(I,J,{\mathbb{G}})$. Given $I, J$ homogeneous ideals in $\mathbb G$, computes $\Delta_J(I)$.
\item $\textsc{GetElement}(L,j)$. It returns the $j$-th element of the list $L$.
\end{itemize}

\begin{algorithm}[H]
\caption{Algorithm computing $g\in \mathrm{PGL}$ and $J\in \mathcal B_{p(t)}$ such that $\mathcal H_J$ is an open subset of the Hilbert scheme containing the point  $\Proj(S/I^g)$}
\label{algo1}
\begin{algorithmic}[1]
\STATE $\textsc{BorelOpenSet}(I)$
\REQUIRE $I$ homogeneous ideal in $K[x_0,\ldots,x_n]$.
\ENSURE the pair $(g,J)$, with $g\in \mathrm{PGL}$ and $J$ a Borel ideal defining an open subset $\mathcal{H}_J$ of the Hilbert scheme containing the point corresponding to $\Proj(S/I^g)$. 
\STATE $p(t)\leftarrow$ Hilbert polynomial of $S/I$;
\STATE $(\tilde I, \mathbb{G}) \leftarrow \textsc{SuitableGrasmannian}(I)$;
\STATE $\mathcal B_{p(t)}\leftarrow\textsc{BorelIdeals}(n,p(t))$;
\STATE $\textsf{openSetFound} \leftarrow \FALSE$;
\STATE $\textsf{output} \leftarrow \emptyset$;
\WHILE{\NOT \textsf{openSetFound}}
\STATE $(g,\tilde I^g)\leftarrow\textsc{QGenericChangeOfCoordinates}(\tilde I)$;
\FOR{$i = 1,\ldots,\vert\mathcal{B}_{p(t)}\vert$ \AND \NOT \textsf{openSetFound}}
\STATE $J \leftarrow \textsc{GetElement}(\mathcal B_{p(t)},i)$;
\IF{$\textsc{Pl\"uckerCoordinate}(\tilde I^g,J,\mathbb{G})\neq 0$}
\STATE  $\textsf{openSetFound} \leftarrow \TRUE$;
\STATE $\textsf{output} \leftarrow (g,J)$;
\ENDIF
\ENDFOR
\ENDWHILE
\RETURN \textsf{output};
\end{algorithmic}
\end{algorithm}

\begin{example}\label{J1}
We now execute Algorithm $\textsc{BorelOpenSet}$ on the ideal $I=(x_2^2,x_1^2)\subset S=K[x_0,x_1,x_2]$ of Example \ref{diffUH}. The Hilbert polynomial of $S/I$ is $p(t)=4$ and the Gotzmann number is $r=4$. We now find a generic change of coordinates $g\in \mathrm{PGL}$ and ideal $J\in \mathcal B_{p(t)}$ such that $(I_{4})^g\in \mathcal H_J$.
There are only two  Borel ideals $J$ in $\mathcal B_4$:
\[
J_1 = (x_2^2,x_2x_1,x_1^3)_{\geq 4}, \quad J_2 = (x_2,x_1^4)_{\geq 4}
\]
After a generic linear change of coordinates $g$ (for instance the one used in Example \ref{diffUH}), we compute the coordinate $\Delta_{J_1}((I_4)^g)$ by writing down the matrix of the coefficients of the 14 homogeneous forms of degree 4 generating $(I_4)^g$. The minor corrresponding to the monomial basis of $J_1$ is non-zero, so $(I_4)^g$ belongs to $\mathcal H_{J_1}$. For $g$ as in Example \ref{diffUH}, $(I_4)^g$ is generated by the $J_1$-marked basis:
\[
G=\left\{x_2^4, x_2^3x_1, x_2^3x_0, x_2^2x_1^2, x_2^2x_1x_0, x_2^2x_0^2, x_1^4, x_2x_1^3,  x_1^3x_0, x_2x_1^2x_0, x_2x_1x_0^2+\frac{1}{2} x_1^2x_0^2\right\}.
\]
\end{example}

 Theorem \ref{fuori} does not extend to monomial ideals $J\in \mathcal Q\setminus \mathcal B$, as the following example shows.
 \begin{example}\label{open0dim}
  Let us consider the constant Hilbert polynomial $p(t)=2$ and the Hilbert scheme $\hilb^2_2$ parameterizing 0-dimensional subschemes in $\PP^2$ of degree 2: in this case $r=2$ and $s=4$. The monomial ideal $J=(x_0^2,x_1^2,x_2^2,x_0x_1)$ is generated by $s$ monomials of degree 2, it is not Borel and  obviously does not belong to $\hilb^2_2$ since it is a primary ideal over the irrelevant maximal ideal $(x_0,x_1,x_2)$. Nevertheless, $\mathcal{H}_J$ is not empy, as it contains   the ideal $(x_0^2-x_0x_2, x_1^2-x_1x_2,x_2^2-x_0x_2-x_1x_2,x_0x_1)$ corresponding to the set of points $\{[1:0:1], [0:1:1]\}$ and, more generally, all the ideals corresponding to  pairs of  distinct points outside the line $x_2=0$ and not on the same line through $[0:0:1]$.
  \end{example}

Theorem \ref{fuori} not only allows us to define a cover of $\hilbp$ from the Borel ideals of $\mathcal B_{p(t)}$, but also gives an interesting role to the ideals belonging to $\mathcal B\setminus \mathcal B_{p(t)}$: indeed, $\hilbp$ lies on every hypersurface of $\mathbb G$ defined by the vanishing of the \Pl\ coordinate  corresponding to such a Borel ideal (and on every hypersurface obtained from these by the action of $\mathrm{PGL}$). 
\begin{lemma}\label{crescitagrado}
For every distinct $J,J'\in \mathcal Q$,  $\mathcal U_J\setminus \mathcal U_{J'}$ is a hypersurface of $\mathcal U_J\simeq \mathbb A^D$ of degree $\vert B_J\setminus B_{J'}\vert$.
\end{lemma}
\begin{proof}
We simply need to observe that the equation of $\mathcal U_J\setminus  \mathcal U_{J'}$ in $\mathcal U_J$ is defined by $\Delta_{J'}/\Delta_J$, whose degree is exactly $\vert B_J\setminus B_{J'}\vert$ by Lemma \ref{lem:pluckerlocali}.
\end{proof}

\begin{corollary}
In the previous setting, we have the set-theoretical inclusion:
\begin{equation*}
\hilbp \subseteq \bigcap_{\begin{subarray}{c}g\in \PGL \\ J\in \mathcal{B}\setminus \mathcal B_{p(t)}\end{subarray}}\Pi^g_J
\end{equation*}
where $\Pi_J$ is the hypersurface of $\mathbb G$ given by $\Delta_J=0$.
\end{corollary}
\begin{proof}
For every $J\in \mathcal B\setminus \mathcal B_{p(t)}$,  $\mathcal H_J$ is empty  by Theorem  \ref{fuori}. This means that the hypersurface $\Pi_J$ contains $\hilbp$. Furthermore, since $\left(\hilbp\right)^g=\hilbp$ for every $g\in \mathrm{PGL}$, we also have $\hilbp\subseteq \Pi_J^g$ for every $J\in \mathcal B\setminus \mathcal B_{p(t)}$ and every $g\in \mathrm{PGL}$.
\end{proof}

If the Hilbert polynomial $p(t)$  is the constant $r$,   then  every Borel ideal $J\in \mathbb{G}$ belongs to $\hilb^n_r$ i.e. $\mathcal{B}=\mathcal B_{p(t)}\subset \hilbp$ (see \cite[Theorem 3.13]{CLMR}).  Then in the 0-dimensional case the family of hypersurfaces $\Pi_J$ considered in Theorem \ref{fuori} is in fact empty. 
If $\deg p(t) = d\geq 1$, $\mathcal{B}\setminus \mathcal B_{p(t)}$  in general is not empty and its elements define subschemes of $\PP^n$ of dimension equal to or lower than the one of the subschemes parametrized by $\hilbp$.
Indeed, if $I\in \mathcal B$  and $S/I$ has Hilbert polynomial $\widetilde{p}(t)\neq p(t)$, then
by Macaulay's lower bound and Gotzmann's persistence, $\dim_K (I_t)>q(t)$ for $t\geq r$. Hence for $t\gg0$, $\widetilde q(t)>q(t)$ and $\widetilde p(t)< p(t)$. Therefore $\deg \widetilde p(t)\leq d$.

\begin{example}\label{nonBorel}
  Let us consider the Hilbert polynomial $p(t)=3t$ in $\PP^3$. The closed points of  $\hilb^3_{3t}$ correspond to curves in $\PP^3$ of degree 3 and arithmetic genus 1, hence it contains all the smooth plane elliptic curves and also some singular or reducible or non-reduced curves.  The  Gotzmann number of $p(t)$ is $r=3$, 
  then $s=q(3)=11$.
  We embed $\hilbp$ in the Grassmannian $\mathbb{G}(11,K[x]_3)$. 
   The only Borel ideal in $\hilb^3_{3t}$ is the $\mathtt{Lex}$-segment ideal:  $$J_{\mathtt{Lex}}=\left(x_3^{3},x_3^{2}x_2,x_3x_2^{2},x_2^{3},x_3^{2}x_1,x_3x_2x_1,x_3x_1^{2},x_3^{2}x_0,x_3x_2x_0,x_3x_1x_0,x_3x_0^{2}\right).$$ 
 The Borel cover of $\hilb^3_{3t}$ is then $\cup_{g\in \mathrm{PGL}} \mathcal H^g_{J_{\mathtt Lex}}$,  with $\mathcal H_{J_{\mathtt Lex}}:=\mathcal{U}_{J_{\mathtt{Lex}}}\cap \hilb^3_{3t}$. 

 We can compute the Borel ideals $J_i$ in $\mathbb{G}(11,K[x]_3)$ that do not belong to $\hilb^3_{3t}$ and the Hilbert polynomial of $S/J_i$:\\
\indent $J_1=(x_2^3,x_3^2,x_2^2x_1, x_3x_1,x_3x_2)_{\geq 3}$, {Hilbert polynomial of } $S/J_1$: $2t+3$;\\
\indent $J_2= (x_2^2,x_3^2,x_3x_1^2,x_3x_2)_{\geq 3}$, {Hilbert polynomial of } $S/J_2$: $2t+3$;\\
\indent $J_3= (x_2^3,x_3^2,x_2x_1^2,x_2^2x_1,x_3x_1^2,x_3x_2)_{\geq 3}$, {Hilbert polynomial of }$S/J_3$: $t+6$;\\
\indent $J_4= (x_1^3,x_2^3,x_3^2,x_2x_1^2,x_2^2x_1,x_3x_1^2,x_3x_2^2,x_3x_2x_1)_{\geq 3}$, {Hilbert polynomial of }$S/J_4$: $9$.\\
  Then $\hilb^3_{3t}$ as a  subscheme of  $\mathbb{G}(11,K[x]_3)$  is contained set-theore\-ti\-cal\-ly in the hypersurfaces $\Pi_{J_i}$ given by $\Delta_{J_i}=0$, $i=1,\dots,4$ (and in every hypersurface obtained from these by the action of $\mathrm{PGL}$). 
  \end{example}


 \section{\texorpdfstring{Equations defining $\mathcal{H}_J$}{Equations defining H\_J}}\label{quarta}

The aim of the present section is to find suitable affine subspaces of  $\mathbb A^D\simeq\mathcal U_J$ in which $\mathcal H_J$ can be isomorphically projected, and furthermore to study in which cases we can control the number and the degree of a set of generators of the ideals defining $\mathcal H_J$ as a closed subscheme of such affine subspaces.  

The results we obtain are very similar to the ones for Gr\"obner\ strata showed by \cite{LR}, but are more general  (see \cite[Example 6.2]{BCLR}). Here the isomorphism of Theorem \ref{isomHM} is crucial, because it allows the use of the techniques  presented in \cite{BCLR} for marked schemes over a $m$-truncation $J_{\geq m}$ of a Borel ideal $J$.

We now consider  a Borel ideal $J$ belonging to $\mathcal B_{p(t)}$: as pointed out in Remark \ref{ossTronc}, in our setting $J$ is an $r$-truncation ideal, in other words $J={J^\sat}_{\geq r}$. We denote by $r'$  the regularity of  $J^\sat$, $r'\leq r$, and by $\rho$ the maximal degree of a monomial in $B_{J^\sat}$ divisible by $x_1$; if there are no such monomials in $B_{J^\sat}$, we set $\rho:=0$. \\
We consider  the  $m$-truncation Borel ideal ${J^\sat}_{\geq m}$. Finally, we will denote by $\phi_{J,r}$ the embedding $\mathcal H_J\hookrightarrow \mathbb A^{p(r)q(r)}$ given by Theorem \ref{isomHM} and by $\mathfrak{A}_J\subset K[C]$ the ideal defining  $\mathcal{H}_J$ as a subscheme of $\mathbb{A}^{p(r)q(r)}$.

\begin{theorem}\label{isomH}
In the above setting, the following statements hold:
\begin{enumerate}[(i)]
\item if $m\geq r$, then  $\Mf({J^\sat}_{\geq m})\simeq \mathcal H_J$;
 \item if $m<r$, then  $\Mf({J^\sat}_{\geq m})$ is a closed subscheme of $\mathcal H_J$, possibly equal. If we consider the embedding  $\phi_{J,r}(\mathcal H_J) \subset \mathbb A^{p(r)q(r)}$, then $\Mf({J^\sat}_{\geq m})$ is cut out by a suitable linear space;
  \item $\mathcal H_J\simeq\Mf({J^\sat}_{\geq m})$ if, and only if, either ${J^\sat}_{\geq m}= J$ or $m\geq \rho-1$.
   \end{enumerate}
 In particular, if $\rho >0$, then  $\rho -1$ is the smallest integer $m$ such that: 
\[ \mathcal H_{J}\simeq\Mf({J^\sat}_{\geq \rho-1})\]
The isomorphism  $\mathcal H_{J}\simeq \Mf({J^\sat}_{\geq r'})$ induces an embedding $\phi_{J,r'}$ of $\mathcal H_{J}$ in an  affine space of dimension   $\vert B_{J^\sat}\vert \cdot p(r')$ and the isomorphism  $\mathcal H_{J}\simeq \Mf({J^\sat}_{\geq \rho-1})$ induces an embedding  $\phi_{J,\rho -1}$ of  $\mathcal H_J$ in an  affine space of dimension 
   \begin{equation}\label{dimmin}
\sum_{x^{\alpha}\in B_{{J^\sat}_{\geq\rho-1}}}\left\vert \cN({J^\sat}_{\geq \rho-1})_{\vert\alpha\vert}\right\vert.
\end{equation}
\end{theorem}
\begin{proof}
Thanks to the isomorphism $\mathcal H_J\simeq \Mf(J)$ of Theorem \ref{isomHM}, the statements are straightforward consequences of \cite[Theorem 5.7]{BCLR}.
\end{proof}

The embeddings  $\phi_{J,\rho -1}$ (or more generally $\phi_{J,m}$ with $\rho-1 \leq m <r'$) of  $\mathcal H_J$ in affine spaces defined in Theorem \ref{isomH} are computationally  advantageous, because in order to compute equations for $\mathcal H_J$ we   deal with a number of variables smaller than $p(r)q(r)$. Furthermore, the ideal defining the affine scheme $\Mf({J^\sat}_{\geq m})$ can be explicitly computed using the algorithm in \cite[Appendix]{BCLR}: this algorithm does not perform any  elimination of variables in order to pass from the ideal defining $\mathcal H_J$ in $\mathbb A^{p(r)q(r)}$ to the ideal defining it as a subscheme of a smaller affine space.

\begin{example}\label{ex:troncatoGrado2}
We consider again the ideal $J_1$ of Example \ref{J1}. We  can compute equations for the ideal defining $\mathcal H_{J_1}$ in the affine space of dimension given by formula  \eqref{dimmin}, namely in this case $12< p(r)q(r)=44$. Observe that $\rho-1$ in this case is 2, hence $\mathcal H_{J_1}$ is isomorphic to $\Mf({J_1}^\sat)$:  by the algorithm presented in \cite[Appendix]{BCLR}, we obtain for $\mathfrak A_{{J_1}^\sat}$  a set of  8 generators  of degree 3 in 12 variables for the ideal.
Starting from the ${J_1}^\sat$-marked set
\[
\begin{split}
& x_2^{2}-C_{1,1}x_1^{2}-C_{1,2}x_2x_0-C_{1,3}x_1x_
0-C_{1,4}x_0^{2}, \\
& x_2x_1-C_{2,1}x_1^{2}-C_{2,2}x_2x_0-C_{2,3}x_1x
_0-C_{2,4}x_0^{2}, \\
&x_1^{3}-C_{3,1}x_1^{2}x_0-C_{3,2}x_2x_0^2-C_
{3,3}x_1x_0^2-C_{3,4}x_0^{3}\\
\end{split}
\]
we obtain the following generators
{\small
\[
\begin{split}
&-C_{2,1}C_{2,2}C_{2,4}-C_{2,2}C_{1,4}+C_{1,2}C_{2,4}+C_{1,1}C_{3,4}-{C_{2,1}}^{2}C_{3,4}-C_{2,3}C_{2,4},\\
&-C_{2,3}C_{2,2}-C_{2,1}C_{2,2}^{2}-C_{2,4}+C_{1,1}C_{3,2}-C_{2,1}^{2}C_{3,2},\\
&C_{1,4}-C_{2,1}C_{2,2}C_{2,3}-C_{2,1}C_{2,4}-C_{2,1}^{2}C_{3,3}-C_{2,3}^{2}+C_{1,1}C_{3,3}+C_{1,2}C_{2,3}-C_{2,2}C_{1,3},\\
&-C_{2,1}^{2}C_{3,1}+C_{1,3}-C_{2,2}C_{1,1}+C_{1,2}C_{2,1}-C_{2,1}^{2}C_{2,2}+C_{1,1}C_{3,1}-2\,C_{2,3}C_{2,1},\\
&C_{2,2}^{2}C_{2,4}-C_{3,3}C_{2,4}+C_{2,1}C_{3,2}C_{2,4}+C_{2,1}C_{2,2}C_{3,4}+C_{2,3}C_{3,4}-C_{3,2}C_{1,4}-C_{3,1}C_{2,2}C_{2,4},\\
&2\,C_{2,1}C_{3,2}C_{2,2}-C_{3,4}-C_{3,3}C_{2,2}-C_{3,1}C_{2,2}^{2}+C_{2,3}C_{3,2}+C_{2,2}^{3}-C_{3,2}C_{1,2},\\
&C_{2,1}C_{3,2}C_{2,3}+C_{2,2}C_{2,4}+C_{2,1}C_{2,2}C_{3,3}-C_{3,1}C_{2,2}C_{2,3}-C_{3,1}C_{2,4}+C_{2,1}C_{3,4}-C_{3,2}C_{1,3}+C_{2,2}^{2}C_{2,3},\\
&C_{2,1}^{2}C_{3,2}-C_{1,1}C_{3,2}+C_{2,1}C_{2,2}^{2}+C_{2,4}+C_{2,3}C_{2,2}.
\end{split}
\]
}

{By Theorem \ref{isomH}, we have embedded $\mathcal H_{J_1}\simeq\Mf({J_1}^\sat)$ in $\mathbb A^{12}$: however, this is not the smallest affine space in which $\Mf({J_1}^\sat)$ can be embedded. Indeed, we can see that from the above set of generators of the ideal $\mathfrak A_{{J_1}^\sat}$, we can easily eliminate the variables $C_{1,3}$, $C_{1,4}$, $C_{2,4}$ and $C_{3,4}$ (each of them appears linearly in a different generator); eliminating them, we obtain the ideal $(0)$, hence $\Mf({J_1}^\sat)$ is isomorphic to $\mathbb A^8$ and this is clearly the smallest affine space in which  $\mathcal H_{J_1}\simeq\Mf({J_1}^\sat)$ can be embedded. Furthermore, since $\mathcal H_{J_1}\simeq \mathbb A^8$, $J_1$ is a smooth point of $\hilb^4_2$.}

\end{example}

From the computational point of view, it would be useful and interesting to bound the number and degree of a set of generators of the ideal defining                     $\mathcal H_J$ given by the embeddings $\phi_{J,m}$ of Theorem \ref{isomH}.

For what concerns a bound on the degree, we might try to derive it, at least for  $\phi_{J,r}$, from other bounds in literature on the degree of the generators for an  ideal $\Gamma$ in $K[\Delta]$ defining the Hilbert scheme $\hilbp$ in $\mathbb G$. 
In \cite{IarroKanev} it is proved that there is such an ideal $\Gamma$ generated in degree $q(r+1)+1$. In \cite{B} it is conjectured, and later on in \cite{HaimSturm} it is proved that, there is such a  $\Gamma$ generated in the far lower degree $n+1$. More recently,  in the case of Hilbert schemes of points, that is $d=0$,  in \cite{ABM} the authors improved the previous ones, showing that there is a $\Gamma$ generated in degree $2$. Finally, in the recent paper \cite{BLMR} the authors proved that the Hilbert scheme can be defined by an ideal $\Gamma$  generated in degree $\leq d+2$. 

However, when we restrict  $\Delta_{J'}$, for every $J'\in \mathcal Q$,   to $\mathbb A^{p(r)q(r)}\simeq\mathcal U_J$, $J\in \mathcal Q$, $\Delta_{J'}$ becomes a polynomial of degree $\leq q(r)$ in $K[C]$.
Thus, restricting anyone of the above ideals $\Gamma$ of $K[\Delta]$ defining $\hilbp$ to $K[C]$, the previous bounds are multiplied by a factor $q(r)$.

We can also try to get a bound for the number  of generators for the ideal $\mathfrak A_J$ defining $\mathcal H_J$ as a subscheme of $\mathbb A^{p(r)q(r)}$, considering the setting of \cite[Proposition C.30]{IarroKanev}: it is enough count the number of $(s'+1)\times(s'+1)$ minors of the coefficient matrix of $(\mathfrak I_J)_{r+1}$, $\id I_J$ generated by a $J$-marked set as in \eqref{JbaseC}. 
 Unfortunately, this number is ${\displaystyle\binom{(n+1)s}{s'+1} \cdot \binom{N(r+1)}{s'+1}}$, where $N(r+1)={\displaystyle\binom{n+r+1}{n}}$.

We now show that if we consider an embedding $\phi_{J,m}$ for $m\geq r'=\reg(J^\sat)$,  we can  get a set of generators for the ideal defining $\mathcal H_J$ in $\mathbb A^{p(m)q(m)}$ such that its cardinality is bounded by a far smaller number  only depending on $p(t)$, $n$ and $m$; moreover,  we get that the degrees of such generators are $\leq d+2$ for every $m\geq r'$ and for every  $n$.

\begin{theorem}\label{gradod+2}
If $J \in \mathcal B_{p(t)}$, for every $m\geq r'$, $ \mathcal{H}_{J}$  is isomorphic to a closed subscheme of  $\mathbb{A}^{p(m)q(m)}$, defined by at most $\left(q(m)(n+1)-q(m+1)\right)\cdot p(m+1)$ polynomials of degree $\leqslant d+2$.
\end{theorem}
\begin{proof}
Using Theorem \ref{isomHM}, $\mathcal H_{J}=\Mf({J})$. Furthermore, thanks to Theorem \ref{isomH} , we have that $\mathcal H_{J}\simeq\Mf({J^\sat}_{\geq m})$ for every $m\geq r'$. 

We next show how to compute a specific set of generators for the ideal $\mathfrak A_{{J^\sat}_{\geq m}}$ defining the scheme structure of $\Mf({J^\sat}_{\geq m})$. 
 We can obtain a set of generators for  this ideal  using a special procedure of reduction. Given the ${J^\sat}_{\geq m}$-marked set $\mathcal G^{(m)}$ (similar to the one in \eqref{JbaseC} of Definition \ref{defin_A}), we consider $V_{m+1} = \{x_i F_\alpha\ \vert\ F_\alpha \in \mathcal G^{(m)},\ x_i \leq \min x^\alpha \}$ and we now use the procedure $\xrightarrow{V_{m+1}}$ of \cite[Definition 3.2, Proposition 3.6]{CR} to reduce $S$-polynomials of elements in $\mathcal G^{(m)}$. We use a  Buchberger-like criterion on $S$-polynomials analogous to the one for Gr\"obner\ bases in order to obtain a set of generators for the ideal $\mathfrak A_{{J^\sat}_{\geq m}}$, by \cite[Theorem 3.12, Theorem 4.1]{CR}. 

 It is enough to consider $S$-polynomials corresponding to a basis of the syzygies of the ideal ${J^\sat}_{\geq m}$ (see \cite[Remark 3.16]{CR}). As $m \geq r' = \reg(J^{\sat})$, there exists such a basis given by pairs of variables (for instance the one in \cite{EK}). Thus, we consider $S$-polynomials of the type $x_iF_\alpha -x_jF_{\alpha'}$ with $x_i x^\alpha=x_j x^{\alpha'}$. If $x_ix^{\gamma}$  is a monomial of ${J^\sat}_{\geq m+1}$ that appears in $x_iF_\alpha -x_jF_{\alpha'}$, then  $x^\gamma \in \cN(J)_m$. By  Lemma \ref{potenze} (\ref{potenze_ii}), we can show that  $x_ix^{\gamma}$ is equal to $x_hx^\beta$, for $x_h =\min(x_ix^{\gamma})<x_i$ and $x^\beta \in B_{{J^\sat}_{\geq m}}$.  Then we can reduce $x_ix^{\gamma}=x_h x^\beta$ by rewriting it as $x_h (x^\beta-F_\beta)$ \cite[Theorem 3.12]{CR}. If some monomial of $x_hx^\beta-x_hF_\beta$ belongs to ${J^\sat}_{\geq m}$, then again we can reduce it using some polynomial $x_{h'}F_{\beta'}$ with $x_{h'}<x_h$.
  
At every step of reduction  a monomial is replaced by a sum of other monomials multiplied by one of the variables $C$. Thus, at every step of reduction the  degree of the coefficients in $K[C]$  directly involved  increases  by 1.  If $x^{\eta_0},x^{\eta_1} \dots, x^{\eta_s}$ is a sequence of monomials of degree $m+1$ in ${J^\sat}_{\geq m}$ such that $x^{\eta_{i+1}}$ appears in the reduction of $x^{\eta_{i}}$, then $\min(x^{\eta_{i+1}}) <\min(x^{\eta_{i}})$.  Since $m \geq r' = \reg(J^{\sat})$, by Lemma \ref{potenze}, (\ref{potenze_iii}) the minimal variable of any monomial in $\cN(J^\sat)_{m+1}$ is smaller than $x_{d+1}$, so the   length of any such  chain   is at most $d+1$. Thus, in the complete reduction of an $S$-polynomial by $\xrightarrow{\ V_{m+1}\ }$ (which is a polynomial whose support is contained in $\cN(J^\sat)_{m+1}$), the final degree of the coefficients in $K[C]$ of each monomial $x^\gamma$  is at most $1+1\cdot(d+1)=d+2$. 

As said before,  it is sufficient to reduce the $S$-polynomials corresponding to a special basis of the syzygies of ${J^\sat}_{\geq m+1}$. If we consider the \emph{Eliahou-Kervaire $S$-polynomials} \cite[Definition 2.12]{BCLR}, we obtain $\left(q(m)(n+1)-q(m+1)\right)$ polynomials. For each of them, using the polynomial reduction process just described, we obtain at most $p(m+1)$ generators for $\mathfrak A_{J^\sat_{\geq m}}$, since the monomials $x^\gamma$ of the complete reduction of a polynomial by $\xrightarrow{\ V_{m+1}\ }$ are contained in $\cN(J^\sat)_{m+1}$ and $m+1 > r'$.
\end{proof}

\begin{example}\label{numerominori}
We consider again $p(t)=4$ in $\PP^2$ as in Example \ref{J1}. In this case we have $r=4$,  $s'=q(5)=17$, $N(5)=21$. The number of $18\times 18$ minors of the coefficient matrix of $(\mathfrak I_J)_{r+1}$ (which has 33 rows and 21 columns), for every $J\in \mathcal Q$, is  $1379420565600$ and their degree is $\leq 18$.
If we consider $J_1$ as in Examples \ref{J1} and \ref{ex:troncatoGrado2}, and use Theorem \ref{gradod+2},  we obtain that $\mathcal H_{J_1}$ is isomorphic to  a subscheme of $\mathbb A^{p(3)q(3)}=\mathbb A^{24}$ defined by at most 28 equations of degree $\leq 2$.
\end{example}

Theorem \ref{gradod+2} allows us to explicitly compute $\phi_{J,m}(\mathcal H_J)$ with a good compromise between the number of involved variables and the degree of the defining equations. This is the reason why in Algorithm $\textsc{BorelOpenSet}(I, p(t),n)$ we sort the list $\mathcal B_{p(t)}$ according to increasing regularity of the saturated ideals. In fact, given $I\in \hilbp$ and a generic $g\in \mathrm{PGL}$, the algorithm returns $J\in \mathcal B_{p(t)}$ such that $I^g\in \mathcal H_J$ with minimum $\reg(J^\sat)$. This turns out to be computationally convenient: in fact, $\mathcal H_J$, among those containing $I^g$,  can be   embedded by equations of degree $\leq d+2$ in an affine space of smallest dimension. Furthermore, by Theorems \ref{isomH} and \ref{gradod+2}, we can improve Algorithm $\textsc{BorelOpenSet}(I, p(t),n)$ by replacing the computation of $\Delta_J(I)$ with an analogous test concerning $\Mf({J^\sat}_{\geq \rho-1})$ or $\Mf({J^\sat}_{\geq \reg(J^\sat)})$. This would be very useful in a practical implementation of this algorithm since it allows to deal with coefficient matrices of a significantly smaller size.

Furthermore, we are able to study the glueing of any pair of subsets of the Borel cover, by choosing a convenient ambient space among those given by the embeddings studied in Theorems \ref{isomH} and \ref{gradod+2}.
\begin{corollary}\label{smallemb} 
If $\phi_{J_i,r}: \mathcal H_{J_i} \rightarrow \mathbb{A}^{p(r)q(r)} $ are the embeddings  for the open subsets corresponding to two  Borel ideals   $J_1$ and $J_2$  belonging to $\hilbp$, then:  
	\[
	\phi_{J_1,r}(\mathcal H_1\cap \mathcal H_2) =\phi_{J_1,r}(\mathcal H_1)\setminus F_1 , \quad \phi_{J_2,r}(\mathcal H_1\cap \mathcal H_2) =\phi_{J_2,r}(\mathcal H_2)\setminus F_2,
	\]
	 where $F_1$ e $F_2$ are hypersurfaces  in $\mathbb{A}^{p(r)q(r)}$ of the same degree $\vert B_{J_1} \setminus B_{J_2}\vert$.\\
  If   $\overline{m}\geq \max \{ \reg({J_1}^\sat), \reg({J_2}^\sat)\}$
	 the same statement   holds for  $\phi_{J_i,\overline{m}}\colon \mathcal H_{J_i} \hookrightarrow \mathbb{A}^{p(\overline{m})q(\overline{m})} $.
\end{corollary}
\begin{proof} 
 $F_1$ is defined by the equation of $\frac{\Delta_{J_2}}{\Delta_{J_1}}$ in $\mathbb A^{p(r)q(r)}$, hence its degree corresponds to $\vert B_{J_1} \setminus B_{J_2}\vert$, by Lemma \ref{lem:pluckerlocali}. The statement follows observing that since $J_1$ and $J_2$ both belong to $\hilbp$,  $\vert B_{J_1} \setminus B_{J_2}\vert =\vert B_{J_2} \setminus B_{J_1}\vert$.
 
The last statement is a straightforward consequence of Theorem \ref{gradod+2}.
\end{proof}

\section{Examples}\label{secexamples}

We now present a few examples to illustrate the twofold interest of the above results, mainly Theorems \ref{isomHM}, \ref{fuori},  \ref{isomH}, and \ref{gradod+2}. On the one hand,   we can  perform explicit computations (by the algorithm in \cite[Appendix]{BCLR}) to embed an open subset of $\hilbp$ in a suitable affine space. On the other hand,  we can also use the results of the present paper to investigate some properties and features of $\hilbp$.

\begin{example}
We consider the Hilbert scheme of $\mu$ points in $\mathbb P^n$, $\hilb_{\mu}^n$. For any monomial ideal $J$ , we have that  the open subset of the Grassmannian $\mathcal U_J$ is isomorphic to $\Af^{p(r)q(r)}$, where $r=\mu$, $p(r)q(r)=\mu\left(\binom{n+\mu}{n}-\mu\right)$.
We also consider the monomial saturated ideal ${L}=(x_n,\dots,x_2,x_1^{\mu})$, which is a ${\mathtt{Lex}}$-segment, its regularity is $r=\mu$ and ${L}_{\geq \mu}$ is in $\hilb_{\mu}^n$.  Obviously, the open subset $\mathcal H_{L_{\geq \mu}}$, the ${\mathtt{Lex}}$-component, contains all the subschemes of $\mathbb P^n$ of $\mu$ distinct points, thus it has dimension $\geq n\mu$. 
By Theorem  \ref{isomH}, $\mathcal H_{L_{\geq \mu}}$ is a closed subscheme of $\mathbb A^{n\mu}$, where $n\mu=\vert B_{L}\vert \cdot p(r)$. Therefore, $\mathcal H_{L_{\geq \mu}}\simeq \Af^{n\mu}$.
\end{example}

\begin{example}
We can now easily study some features of $\hilb_{3t}^3$ that we have already investigated in Example \ref{nonBorel}. 
The open Borel cover of $\hilb_{3t}^3$  is made up of the open subsets $\mathcal H^g_{J_{\mathtt{Lex}}}$, $g\in \mathrm{PGL}$ and $J_{\mathtt{Lex}}=(x_3,x_2^3)_{\geq 3}$, as already pointed out in Example \ref{nonBorel}, using Theorem \ref{fuori} and  Theorem \ref{isomHM}. Since no monomial in the basis of $(x_3,x_2^3)$ is divisible by $x_1$, using Theorem  \ref{isomH}, we have that $\Mf((x_3,x_2^3))\simeq \mathcal H_{J_{\mathtt{Lex}}}\hookrightarrow \Af^{12}$. By explicit computations, we obtain that the embedding is actually an isomorphism, that is $\mathcal H_{J_{\mathtt{Lex}}}\simeq \Af^{12}$.

Indeed, every ideal $I$ in $\Mf((x_3,x_2^3))$ is generated by a linear form $x_3+L$ and a cubic form $Q$, $L,Q \in K[x_0,x_1,x_2]$, hence it depends on 12 free parameters.
\end{example}

\begin{example}
Consider $n=2$ and $p(t)=7$; in this case, we consider the set  $\mathcal B_{p(t)}$  in $K[x_0,x_1,x_2]$:
\begin{itemize}
\item $J_1 = (x_2,x_1^7)_{\geq 7}$;
\item $J_2 = (x_2^2,x_2x_1,x_1^6)_{\geq 7}$;
\item $J_3 = (x_2^2,x_2x_1^2,x_1^5)_{\geq 7}$;
\item $J_4 = (x_2^2,x_2x_1^3,x_1^4)_{\geq 7}$;
\item $J_5 = (x_2^3,x_2^2x_1,x_2x_1^2,x_1^4)_{\geq 7}$.
\end{itemize}

$\hilbp$ in this case can be covered by the open subsets $\mathcal H_{J_i}$, $i=1,\dots,5$ (up to changes of coordinates, as shown in Theorem \ref{fuori}), while the cover obtained  from the complete list of monomial ideals in  $\mathbb{G}(29,36)$ consists of 8347680 open subsets, each one isomorphic to   $\mathbb{A}^{203}$.

If we want to compute the equations defining $\mathcal H_{J_i}\simeq \Mf(J_i)$ as an affine scheme, we can use the techniques developed in \cite{BCLR}. Hence we can choose whether we are interested in embedding $\mathcal H_{J_i}$ in the smallest affine space (by considering the isomorphism with $\Mf({J_i^\sat}_{\geq m})$, for some $m$ greater than or equal to $\rho_i-1$, $\rho_i$ as in Theorem  \ref{isomH}) or wheter we also want to keep control on the degree of the equations. Indeed in a  larger affine space, $\mathcal H_{\mathfrak j_i}$ is defined by equations of degree $\leq 2$ (see Theorem \ref{gradod+2}).

It is interesting to point out  that $J_4$ does not fullfill the hypotheses of results  by \cite[Section 6]{LR}, so these results for Gr\"obner Strata do not apply to  $\mathcal H_{J_4}$, while the techniques of \cite{BCLR} do.  Similarly, in \cite[Appendix]{CR}, the authors exhibit an open subset  $\mathcal H_J$ of $\hilb_8^2$, showing  by explicit computations that it is not a Gr\"obner\ stratum.

Finally, the points of $\hilb_7^2$ are the hyperplane sections of the curves of $\mathbb P^3$ parametrized by the Hilbert polynomial $7t-5$. 
This Hilbert scheme is investigated in detail by \cite{GFR} using the techniques of \cite{BCLR} and of the present paper. The  Borel monomial ideals in $K[x_0,x_1,x_2,x_3]$ giving the Borel cover of the Hilbert scheme $\hilb^3_{7t-5}$ is 112. Among these, some do not match the hypothesis of \cite{LR}, while all of them can be handled with the strategies presented in the present paper and in \cite{BCLR}.
\end{example}

For other examples, including the explicit computation of equations for marked schemes, we refer to \cite[Example 6.3]{BCLR}.

\section*{Acknowledgments} 
The authors were supported by PRIN 2008 (Geometria della variet\`a algebriche e dei loro spazi di moduli) co-financed by MIUR.
The third author thanks the American Institute of Mathematics for financial support and the organizers and participants of the Meeting \lq\lq Components of Hilbert Schemes, Palo Alto 2010\rq\rq\ for useful discussions.

\end{document}